\documentclass[11pt]{article}

\usepackage{times}
\usepackage{alltt}
\usepackage{amssymb}

\usepackage{amssymb,amsmath,latexsym}    % Use special AMS-math environments

\usepackage{array}      % More control ov

\newtheorem{thm}{Theorem}

\newtheorem{lem}{Lemma}

\newtheorem{prop}{Proposition}

\newtheorem{defn}{Definition}
\newtheorem{exmp}{Example}
\newtheorem{rem}{Remark}

\title{Orders of vanishing of zeros of characteristic $p$ zeta function}
\author{Victor Bautista-Ancona, Javier Diaz-Vargas, and Jos\'e Luis Maldonado-Baz\'an}
\date{Facultad de Matem\'{a}ticas.
       Universidad Aut\'{o}noma de Yucat\'{a}n, M\'exico.}

\begin{document}

\maketitle

\section*{Abstract}

Orders of vanishing of zeros of zeta functions have much arithmetic information encoded in them. For the absolute zeta function, Dinesh Thakur gave sufficient conditions for the order of vanishing of its zeros when the finite field has two elements. Such conditions consider only principal ideals. This result was generalized by Thakur and Diaz-Vargas. Now the conditions involve not only the principal ideals but all the classes of ideals, still in the field of two elements. In this work, we generalize these results to arbitrary finite fields, using similar proofs of Thakur and Diaz-Vargas.

\section {Introduction}

One of the most important topics in the study of zeta functions is the order of vanishing of its zeros. Some results have been found for the characteristic $p$ zeta function and the ``trivial'' zeros that we analyze. In \cite{cpzf}, Thakur gave sufficient conditions for a hyperelliptic function field over the finite field $\mathbb{F}_{q}$, $q=2$, to have order of vanishing $2$ at the negative integer $-s$. An interesting phenomenon is that such conditions involve the sum of the digits in the expansion base $2$ of $s$, $l_2\left(s\right)$ (see theorem \ref{hiper}). This result was generalized by Thakur and Diaz-Vargas in \cite{ozcpzf} considering now  all the ideal classes in the definition of the zeta function, and not only the principal ideals. The conditions to have order of vanishing at least $2$, depend again on the decomposition base $2$ of $es$, where $e$ is the exponent of the ideal class group (see theorem \ref{generalization}).

 In \cite{cpzf}, Thakur says succinctly how to deal with the general case $q=p^n$ and arbitrary function fields, when one considers only principal ideals. The conditions in order to have multiplicity $q$ at the zeros, depend now on Weierstrass gaps at $\infty$. We add an extra condition and give proof of theorem \ref{dinesh}, which is a generalization of Thakur's theorem for $q=2$. We analyze also what it means for a function field to have an $r$-gap structure at $\infty$.

 Finally, in theorem \ref{tesismc} we obtain, by joining theorems \ref{dinesh} and \ref{generalization}, our principal result. This is a partial generalization to any value of $q$ and considering all the classes of ideals, of conditions that guarantee to have multiplicity at least $q$ at the trivial zeros, for a very specific class of function fields.

We present also examples of the phenomena described in the theorems, wherever possible. 

\section {Basic definitions}

Let $K$ be a function field of one variable with constant field $\mathbb{F}_{q}$ of characteristic $p$, $\infty $ a place of degree $1$. Then its residue field at $\infty$ is $\mathbb{F}_{q}$ and its completion $K_{\infty }$ at $\infty$ is the Laurent series field $\mathbb{F}_{q} \left( \left( u_{\infty} \right) \right)$, where $u_{\infty}$ is any uniformizer at $\infty$, i.e., a generator of the maximal ideal of the valuation ring at $\infty$.

Let $A$ be the ring of elements of $K$ that have no poles outside $\infty $. Let $\mathbb{C}_\infty $ the completion of the algebraic closure of $K_{\infty }$. We denote by $h$ the class number of $K$, so the class number of $A$ is $h$ too. 

For $x\in K,$ we define the degree of $x$ as $\deg\left( x\right) =-v_{\infty }\left( x\right) .$ By convention, $\deg\left( 0\right) =-\infty $. Similarly, for $\mathcal{I}$ ideal $\in A$ define $\deg \mathcal{I}=\dim_{\mathbb{F}_q}A/\mathcal{I}$.

We define now the relevant zeta function (see \cite{ffa}, p. $156$). 

\begin{eqnarray*}
\zeta \left(s,X\right)&:=&\zeta_{A} \left(s,X\right):=\sum_{d=0}^{\infty}X^{d} \sum_{\deg \: a=d \atop a \: monic} \frac{1}{a^s} \in K\left[ \left[X\right] \right] , \\
\zeta \left(s \right)&:=&\zeta _{A} \left( s\right):=\zeta \left( s,1\right)  \in K_{\infty}.
\end{eqnarray*}

In this work, we will focus on orders of vanishing of trivial zeros of the characteristic $p$ absolute zeta function. We begin by giving the definition of order of vanishing.

\begin{defn}\label{orderv}
The order of vanishing of $\zeta$ at the zero $s$, is defined as the order of vanishing of $\zeta \left( s,X\right)$ at $X=1$.
\end{defn}

\section{Zeta values at negative integers}

For a non negative integer $k=\sum k_{i}q^{i}$, with $0\leq k_{i}<q$, we let $l_{q}\left( k\right)=\sum k_{i}$, that is, $l_{q}\left( k\right) $ is the sum of base $q$ digits of $k$.

\begin{thm} [\cite{cpzf}]
\label{ev}Let $W$ be a $\mathbb{F}_{q}$-vector space of dimension $d$ inside a field $\left( \textrm {or ring}\right) $ $\mathcal{F}$ over $\mathbb{F}_{q}$. Let $f \in \mathcal{F}-W$. If $d>\frac{l_{q}\left( k\right) }{q-1}$, then
\[
\sum_{w\in W}\left( f+w\right) ^{k}=0.
\]
\end{thm}

This proposition is useful for proving the following result, about the ``trivial zeros'' of the absolute zeta function. The theorem is analogue to the classical result, since multiples of $q-1$ are playing the role, in this context, of even integers. 

\begin{thm}[\cite{vzflsmff}]
\label{zetacero}For a positive integer $s$, $\zeta \left( -s\right) \in A$ and 
$\zeta \left( -\left( q-1\right) s\right) =0$.
\end{thm}

The reciprocal is true in the following situation.

\begin{thm}[\cite{vzflsmff}]
\label{equiv}Let $A=\mathbb{F}_{q}[x],$ and let $s$ be a positive integer. Then $\zeta \left( -s\right) =0$ if and only if $s$ is a multiple of $q-1$.
\end{thm}

It is not known for a general $A$ if the values at odd integers are not zero.

\section{Order of vanishing} \label{zeros}

An immediate consequence of theorem $\ref{zetacero}$ is that the negative integers are zeros of the zeta function for $A=\mathbb{F}_2[x]$. In fact more is true; from definition \ref{orderv}, we get the following

\begin{thm}[\cite{zm}]
The negative integers are zeros of the zeta function for $A=\mathbb{F}_{2}\left[ x\right]$, whose order of vanishing is one.
\end{thm}  

\subsection{Considering only principal ideals: case $q=2$} \label{idelprin}

In this and in the following subsections we will give some theorems related with the order of vanishing of the zeta function. We will do it in two ways, first considering only principal ideals and then considering all the ideals (principal or not). Observe that in the first case, the absolute zeta function is a complete zeta function.

We have the following, easy to prove, lemma.

\begin{lem}
\label{lem}
We have 
\[W_{d}=\left\{ a\in A\mid \deg \left( a\right)<d\right\} =\mathcal{L}\left( \left( d-1\right) \infty \right).\]
\end{lem}

Then a simple application of the theorem \ref{ev} and the Riemann-Roch theorem to $W_d$, together with the simplicity of the zeros gives \cite{cpzf} the

\begin{thm}[Thakur]\label{hiper}
If $q=2$ and $K$ is hyperelliptic, then the order of vanishing of $\zeta \left(- s\right)$, $s$ a positive integer,  is $2$ if $l_{2}\left(s\right) \leq g$, where $g$ is the genus of $K$.
\end{thm}

\begin{exmp}

Let $A=\mathbb{F}_2[x,y]$ where $y^{2}+(x^{2}+x+1)y+(x^{2}+x+1)(x^{5}+x^{2}+1)=0.$

\end{exmp}

This is example $26$ of $\cite{qffe2icg}$; we have that $g=3$. So, by theorem $\ref{hiper}$, the order of vanishing is $2$ if $s=(2^l+2^m+2^n), l \geq 0, m\geq 0, n \geq 0$ (possibly with $l=m=n$).$\hfill\Box$
   
Theorem $\ref{hiper}$ is not valid for $q>2$, that is, we can not assure that if $s=k(q-1)>0$, and $-s$ is a zero of $\zeta$ with $l_q(s) \leq g$, then $-s$ is a zero with order of vanishing $2$. Let's see it in the next example. Let $S_{A}(d)$ be the coefficient of $X^{d}$ in $\zeta _{A}(-s,X)$.

\begin{exmp}

Consider $A=\mathbb{F}_3[x,y]$ where $y^2=x(x+1)(x+2)(x^2+1)$.

\end{exmp}

This is example $36$ from $\cite{qffe2icg}$; we have that $g=2$, $\deg\left(x\right)=2$, and $\deg\left(y\right)=5$. By theorem $\ref{zetacero}$, $\zeta\left(-2\right)=0$. Now, $l_3\left(2\right)=2=g$. Then $W_3$ has dimension $2$, and so we have $l_3\left(2\right)>\dim W_3/\left(q-1\right)=1$, and by theorem $\ref{ev}$, $S_{A}\left(d\right)=0$ for $d \geq 3$ Therefore,

\[\zeta_{A}\left( -2,X\right) = \left(1+X\right)\left(1-X\right).\]

Then, the order of vanishing of $-2$ is $1$.\hfill $\Box$

For more background and perspective on zeta functions, see the books by Goss \cite{bsffa}, and Thakur \cite{ffa}.

\subsection{Considering only principal ideals: case $q=p^n$}

Theorem \ref{hiper} gives us sufficient conditions for a hyperelliptic function field to have order of vanishing $2$ when $q=2$. Now, we try to generalize this result to the case $q=p^n$. For that, we need to define what it means for $K$, to have an $r$-gap structure at $\infty$.

\begin{defn}
Let $g$ the genus of $K$ over $\mathbb{F}_{q}$. We say that $K$ has an $r$-gap structure at $\infty$ if there is an $r$ such that $l(iq\infty)=i+1$, for $1\leq i \leq r$, $l \left( (g+r)\infty \right)=r+1,$ and $rq \leq g+r$.
\end{defn}

\begin{exmp}

Let $K$ be a function field of genus $g=5$ over $\mathbb{F}_{3}$. Then $K$ has not a $1$-gap structure at $\infty$.
\end{exmp}

Suppose that $K$ has a $1$-gap structure at $\infty$, i.e., $l\left( 3\infty \right) =2$, and $l\left( 6\infty \right) =2.$ Consider the following diagram:

\begin{center}
\begin{tabular}{ccccccccc}
1 & 2 & 3 & 4 & 5 & 6 & 7 & 8 & 9 \\
1 &   & 2 &   &   & 2 &   &   &
\end{tabular}
\end{center}

\noindent where the first row is filled with the possible gaps (that lies between $1$ and $2g-1=9$) and  the second row with the values of $l\left( i\infty \right).$

Recall that $i$ is a gap if and only if $l\left( \left( i-1\right) \infty \right) =l\left( i\infty \right),$ and there are $g$ gaps. We know that $1$ is always a gap and that $4, 5, 6$ are gaps by assumption.

Then

\begin{center}
\begin{tabular}{ccccccccc}
(1) & 2 & 3 & (4) & (5) & (6) & 7 & 8 & 9 \\
 1  &   & 2 &  2  &  2  &  2  &   &   &
\end{tabular}
\end{center}

Observe that $3$ is a gap since if it were not a gap then an element of degree $6$ could exist, which does not exist because we know that $6$ is a gap. By the same reason, $2$ is a gap, and $1$ is always a gap. Then $K$ has gaps at $1, 2, 3, 4, 5, 6$, which is a contradiction because it has only $g=5$ gaps. So there is no such a function field.$\hfill\Box$

The following example shows that if the function field $K$ has genus $5$ over $\mathbb{F}_{3}$, and $K$ has a $2$-gap structure, then it has an element of degree $3$.

\begin{exmp}
Let K be a function field of genus $g=5$ over $\mathbb{F}_{3}$ that has a $2$-gap structure at $\infty$. Then its gaps are: $1,2,4,5,7$.
\end{exmp}

In order to prove this a $2$-gap structure means that: $l\left( 3\infty \right) =2,l\left( 6\infty \right) =l\left( 7\infty \right) =3.$ Then $7$ is a gap. Consider the following diagram:

\begin{center}
\begin{tabular}{ccccccccc}
(1) & 2 & 3 & 4 & 5 & 6 & (7) & 8 & 9 \\
 1  &   & 2 &   &   & 3 & 3 &   &
\end{tabular}
\end{center}

Now assume that $2$ is not a gap, then $4,6$ and $8$ are not either. Therefore, 

\begin{center}
\begin{tabular}{ccccccccc}
(1) & 2 & 3 & 4 & 5 & 6 & (7) & 8 & 9 \\
 1  &   & 2 & 3 &   & 3 &  3  & 4  &
\end{tabular}
\end{center}

We observe that this is a contradiction, because this force to $l\left( 5\infty \right) =3$, saying that $6$ is a gap, which is absurd. Therefore, $2$ is a gap.

\begin{center}
\begin{tabular}{ccccccccc}
(1) & (2) & 3 & 4 & 5 & 6 & (7) & 8 & 9 \\
 1  &  1  & 2 &   &   & 3 &  3  &   &
\end{tabular}
\end{center}

From the previous diagram, we see that $3$ is not a gap. Suppose now that $4$ is not a gap, then $8$ is not either.

\begin{center}
\begin{tabular}{ccccccccc}
(1) & (2) & 3 & 4 & 5 & 6 & (7) & 8 & 9 \\
 1  &  1  & 2 & 3 &   & 3 &  3  & 4  &
\end{tabular}
\end{center}

This forces $5$ and $6$ to be gaps. This is a contradiction, since there is an element of degree $3$. Then $4$ is a gap.

\begin{center}
\begin{tabular}{ccccccccc}
(1) & (2) & 3 & (4) & 5 & 6 & (7) & 8 & 9 \\
 1  &  1  & 2 &  2  &   & 3 &  3  &   &
\end{tabular}
\end{center}

Lastly, assume that $5$ is not a gap,

\begin{center}
\begin{tabular}{ccccccccc}
(1) & (2) & 3 & (4) & 5 & 6 & (7) & 8 & 9 \\
 1  &  1  & 2 &  2  & 3 & 3 &  3  &   &
\end{tabular}
\end{center}

and we have that $6$ is a gap, which is absurd since there is an element of degree $3$, that again must exists. In conclusion, we have the complete diagram:

\begin{center}
\begin{tabular}{ccccccccc}
(1) & (2) & 3 & (4) & (5) & 6 & (7) & 8 & 9 \\
 1  &  1  & 2 &  2  &  2  & 3 &  3  & 4 & 5
\end{tabular}
\end{center}

and the gaps are: $1,2,4,5,7.$$\hfill\Box$

The following theorem ``generalizes'' the previous example.

\begin{thm}
Let $K$ be of genus $g$ over $\mathbb{F}_{q}$. Suppose that $K$ has a $r$-gap structure at $\infty$ with $r \geq  q-1$. Then $K$ has an element of degree $q$. Moreover, $K$ does not have elements of degree $\neq iq$ between $1$ and $rq$.
\end{thm}

\textbf{Proof: }
Since $K$ has an $r$-gap structure, we have the following diagram:

\begin{center}
\begin{tabular}{ccccccccccccc}
(1) &  2 & \ldots & q-1 & q & q+1 & \ldots & 2q & \ldots & rq & rq+1 & \ldots&  g+r \\
 1  &    &        &     & 2 &     &        & 3  &        & r+1&      &       &  r+1
\end{tabular}
\end{center}

Observe that by definition, $K$ has as gaps $rq+1, \ldots, g+r$ since $rq\leq g+r$ and $l(rq\infty)=r+1$ and $l \left( (g+r)\infty \right)=r+1$. Therefore, we have $g-r\left( q-1\right)$ gaps. Notice also that between $1$ and $rq$ there are exactly $r$ numbers that are not gaps. Suppose that exists an element of degree $1<k\leq q-1$, that is, suppose that $k$ is not a gap. Let $c$ be the greatest number such that $ck\leq rq$. If $c\leq r$, then $rq=c\left( q-1\right) +res\leq r\left( q-1\right) +res$ with $res<q-1.$ Then $rq\leq r\left( q-1\right) +res$ implies that $r\leq res$, this is a contradiction with the fact that $r\geq q-1.$  Therefore, $c>r$ and this implies that we have more than $r$ numbers that are not gaps, this is a contradiction. Therefore, $1, 2, \dots, q-1$ are gaps. We have the following diagram:

\begin{center}
\begin{tabular}{ccccccccccccc}
(1) & (2) & \ldots & (q-1) & q & q+1 & \ldots & 2q & \ldots & rq & (rq+1)& \ldots&  (g+r) \\
 1   &     &        &      &  2  &     &        & 3  &        & r+1&       &       &  (r+1)
\end{tabular}
\end{center}

If we assume that $q$ is a gap, then between $q+1$ and $2q$ must be a non gap, which implies that there are $q-1$ gaps. The same is true for $2q+1$ and $3q$, etc. So, there are $r(q-1)+1$ gaps, a contradiction. Therefore, $q$ is not a gap and from here $2q, 3q, \ldots, rq$ are also no gaps.  \hfill $\Box$

In the previous theorem the hypothesis that $r \geq q-1$ assures us the existence of an element of degree $q$.

The next theorem is essentially due to Dinesh Thakur (compare with \cite {cpzf}).

\begin{thm}
If $q=p^n$ and $K$ is of genus $g$ with an $r$-gap structure at $\infty$ with $r \geq  q-1$, then the order of vanishing of $\zeta (-s)$ at positive integers $s$, $s$ a multiple of $q-1$, is $q$ if $\frac{l_{q}(s)}{q-1} \leq r$.\label{dinesh}
\end{thm}

\textbf{Proof: }
Let $x$ be an element of degree $q$ (it exists by the previous theorem). Recall that $S_{A}(d)$ is the coefficient of $X^{d}$ in $\zeta _{A}(-s,X)$, and that $W_{d}=\mathcal{L} \left( (d-1) \infty \right)$ (lemma \ref{lem}). For $d>\frac{l_{q}(s)}{q-1}+g$, by Riemann-Roch's theorem (being $h=\textrm {dim} \left(H-(d-1)\infty \right)$, $H$ a canonical divisor) we have
\[
\textrm {dim }W_{d}=\textrm {dim } \mathcal{L} \left( (d-1) \infty \right)= l \left( (d-1) \infty \right)=(d-1)+1-g+h\\
=d-g+h > \frac{l_{q}(s)}{q-1}.
\]
Then $S_{A}(d)=0$ for $d>\frac{l_{q}(s)}{q-1}+g$. So, for $d>g+r$, as $r\geqslant \frac{l_{q}(s)}{q-1}$ we have that $S_{A}(d)=0$. Since $A$ does not have elements of degree $\neq iq$ we only consider, in the summation, elements of degree $dq$. From here, 
\[
\zeta _A (-s,X)=\sum_{d=0}^{r} S_{A}(dq) X^{dq}.
\]
Now, as $x$ has degree $q$ in $A$ but one in $\mathbb{F} _{q} [x]$ we have that $S_{A}(dq)=S_{\mathbb{F}_{q}[x]}(d)$ for $0\leq d \leq r$.
As $\frac{l_{q}(s)}{q-1} \leq r$ then $S_{\mathbb{F}_{q}[x]}(d)=0$ for $d>\frac{l_{q}(s)}{q-1}$. From here, we have
\[
\zeta_{A} (-s,X)=\zeta_{\mathbb{F}_{q}[x]}(-s,X^{q}).
\]
We know that the order of vanishing at zeros for $\zeta_{\mathbb{F}_{q}[x]}(-s,X)$ is one and therefore is $q$ for $\zeta_{\mathbb{F}_{q}[x]}(-s,X^{q})$ and from here, we get the result.\hfill $\Box$

An example of the situation described in the theorem is the following, which appears in \cite{cpzf}, and is due to Jos\'e Felipe Voloch:

\begin{exmp}
\label{voloch} Let $q=3$ and consider $x^2(x-1)^3-y^2+yx^3+y^3+y^5=0$ with $\infty$ being $(1,0)$.  
\end{exmp}

We found the next example using the computer program Kant/Kash.

\begin{exmp}
Let $q=3$ and consider $y^5+x^3 y^3+x^2+x+1=0$.
\end{exmp}

In this case, $q=3,g=6$ and we take $r=2$. The gaps are $1,2,4,5,7,8$. Furthermore, $l(3 \infty)=2$ and $l(6 \infty)=l(8 \infty)=3$. Therefore, we have a $2$-gap structure.
If we take $s=3^m+3^n$, then the order of vanishing of the function in $-s$ is $3$. \hfill $\Box$

\subsection{Considering all the ideals: case $q=2$} \label{idelnonprin}

Now, we need to define a zeta function that involves not only the principal ideals, but all the ideals of $A$. We do it as it is suggested in \cite{cpzf}. Let $e$ be the exponent of the ideal class group of $A $. Let $s$ be a multiple of $e$, and define $\mathcal{I}^{s}$ as $a^{s/e}$, where $a$  is the monic generator of $\mathcal{I}^{e}$. Then, we can define the zeta function as follows: for $s$, an integer multiple of $e$,
 
\begin{eqnarray*}
\zeta \left( s,X\right) &:&=\zeta_{A}\left( s,X\right)
:=\sum_{d=0}^{\infty }X^{d}\sum_{\deg \: \mathcal{ I}=d \atop I \: \mathrm{ideal} \: \mathrm{of}  \: A}\mathcal{I}^{-s}, \\
\zeta \left( s\right) &:&=\zeta _{A}\left( s\right) :=\zeta \left( s,1\right).
\end{eqnarray*}

Thakur and Diaz-Vargas, in \cite{ozcpzf}, generalize theorem \ref{hiper}, including now in the calculation of the zeta function all the ideals and not only the principal ones. 

\begin{thm}
\label{generalization}Let $e$ be the exponent of the ideal class group and $q=2$. Let $A=\mathbb{F}_2\left[ x,y\right] $ given by $y^{2}-a\left( x\right) y=b\left( x\right)$ or $y^{2}-y=b\left( x\right) $, where degree of $x=2$ and degree of $y=N$ is an odd number. Assume that $\mathcal{I}_{k}$, $k=1$, $\ldots $, $h-1$ are the integral ideals representing all the non trivial classes of ideals, with $\deg(\mathcal{I}_{k})=d_{k}$ or order $e_{k}$, as an element of the ideal class group. Moreover, assume that $\mathcal{I}_{k}^{e_{k}}=f_{k} $ is an irreducible polynomial, and that it divides to $b\left( x\right) $. If $N>2\mu +e_{k}d_{k}$ for all $k$, then $\zeta \left( -es\right)$, where $s$ is a positive integer, has order of vanishing at least two, if $l_{q}\left(es\right) \leq \mu $.
\end{thm}

The next example satisfies the conditions of the theorem (for details see \cite{ozcpzf}).

\begin{exmp}

\label{h=4,g=3,ex1}Let $A=\mathbb{F}_{2}\left[ x,y\right] /\left( y^{2}+x(x+1)y+x(x+1)(x^5+x^3+x^2+x+1)\right)$.

\end{exmp}

\subsection{Considering all the ideals: case $q=p^n$}

The main result of this work is in the following theorem which is a generalization of theorem \ref{dinesh}, but now considering all the ideals. The proof is a partial ``generalization'' of the proof of theorem \ref{generalization}.

\begin{thm}
\label{tesismc}Let $e$ be the exponent of the ideal class group and $q=p^n$. Let $K=\mathbb{F}_q(x,y)$ given by $y^q-a(x)^{q-1}y=b(x)$ with a $r$-gap structure at $\infty$. Suppose that degree $y=N>q \deg{a(x)}$ is relatively prime to $p$ and let $u=\frac{b(x)}{a(x)^{q}}$. Let $m(x)$ be an irreducible polynomial in $\mathbb{F}_q[x]$ such that, if $u=m(x)^nz$, $z \in \mathbb{F}_q(x)$, $m(x)$ does not divide $z$, and $n<0$, then $n$ is relatively prime to $p$. Assume that $\mathcal{I}_k$, $k=1,2,\ldots,h-1$ are integral ideals representing all the non trivial ideal classes, with $\deg\left(\mathcal{I}_{k}\right)=d_{k}$, and order $e_{k},$ as an element of the ideal class group. Moreover, assume that $\mathcal{I}_k ^{e_k}=f_k$ is irreducible and divides $b(x)$. If $N>q\mu +e_k d_k$ for all $k$, then $\zeta (-es)$, where $s$ is a positive integer, $s$ a multiple of $q-1$, has multiplicity at least $q$, if $\frac{l_{q}(es)}{q-1} \leq \mu$. 
\end{thm}

\textbf{Proof: }
It follows from (\cite{aff}, p. 117) that $x$ has order $q$ and $\infty$ is a place of degree $1$. Furthermore, we have
\[
\zeta \left( -es,X\right) =\sum_{d=0}^{\infty }X^{d}\sum_{\deg \: \mathcal{I}=d} \mathcal{I}^{es}=\sum_{d=0}^{\infty }X^{d}\sum_{\deg \: a=d \atop a \: monic}a^{es}+\sum_{k=1}^{h-1}\sum_{d=0}^{\infty }X^{d}\sum_{\deg \: \mathcal{I}=d \atop \mathcal{I}\sim \mathcal{I}_{k}^{-1}}\mathcal{I}^{es}. 
\]

Since $K$ has a $r$-gap structure, the principal part 
\[
\sum_{d=0}^{\infty }X^{d}\sum_{\deg \: a=d \atop a \: monic}a^{es} 
\]
contributes to the order of vanishing in $q$ by theorem \ref{dinesh}. In fact, it was proved that
\[
\sum_{d=0}^{\infty }X^{d}\sum_{\deg \: a=d \atop a \:  monic}a^{es}=\zeta _{\mathbb{F}_{q}\left[ x\right] }\left( -es,X^{q}\right).
\]
Also, $\mathcal{I}_{k}^{e_k}=f_k$ has deg $=e_kd_k<N$, and so $f_k=f_k\left(x\right)$, that is, $f_k \in \mathbb{F}_q\left[x\right]$. Write $\frac{e}{e_{k}}=s_{k}$. From here, the $k$- term is
\begin{equation*}
\sum_{d=0}^{\infty }X^{d}\sum_{\deg \: \mathcal{I}=d \atop \mathcal{I}\sim \mathcal{I}_{k}^{-1} }\mathcal{I}^{es}=\frac{1}{f_{k}^{ss_{k}}}\sum_{d=0}^{\infty }X^{d}\sum_{\deg \: \mathcal{I}=d \atop \mathcal{I}\sim \mathcal{I}_{k}^{-1} }\left( \mathcal{I}_{k}\mathcal{I}\right) ^{es}.  \
\end{equation*}

We want to show that the vanishing order, for each $k$, is $q$ (hence, the total order of vanishing is at least $q$). Here, $f_{k}^{ss_{k}}$ is independent of $d$ and therefore, we can ignore it for vanishing considerations.

Now, $\mathcal{I}_{k}\mathcal{I}$ are integral of degree $d+d_{k}$ of the form $g_0 (x)+ y g_1(x)+ y^2 g_2 (x)+ \cdots + y^{q-1} g_{q-1}(x)$ if and only if $\frac{g_0 (x)+ y g_1(x)+ y^2 g_2 (x)+ \cdots + y^{q-1} g_{q-1}(x) }{\mathcal{I}_{k}}=\mathcal{I}$ is integral, if and only if, \\ $\frac{\left( g_0 (x)+ y g_1(x)+ y^2 g_2 (x)+ \cdots + y^{q-1} g_{q-1}(x)   \right)^{e_{k}}}{\mathcal{I}_{k}^{e_{k}}}=\frac{\left( g_0 (x)+ y g_1(x)+ y^2 g_2 (x)+ \cdots + y^{q-1} g_{q-1}(x) \right) ^{e_{k}}}{f_{k}}\in A$, i.e.,
 
\begin{equation*}
\left(  g_0 (x)+ y g_1(x)+ y^2 g_2 (x)+ \cdots + y^{q-1} g_{q-1}(x)   \right)^{e_{k}} \equiv 0\pmod{f_{k}}.
\end{equation*}

Then, $y^q-a(x)^{q-1}y=b(x)$ implies that

\[y^{q+j} \equiv a(x)^{q-1}y^{1+j} \pmod{f_{k}}\]

and therefore, 
 
\begin{eqnarray*}
\lefteqn{\left(  g_0 (x)+ y g_1(x)+ y^2 g_2 (x)+ \cdots + y^{q-1} g_{q-1}(x)   \right)^{e_{k}} \equiv} \\ & & g^{e_k}_0 (x)+ y g'_1(x)+ y^2 g'_2 (x)+ \cdots + y^{q-1} g'_{q-1}(x)  \pmod{f_{k}}. 
\end{eqnarray*}

So, if this is congruent to zero $\pmod{f_{k}}$, $f_k \mid g_0^{e_k}$ which implies that $f_k \mid g_0$, since $f_k$ is an irreducible polynomial in $x$. So, we have 
\[
\sum_{ \deg \: \mathcal{I}=d \atop \mathcal{I}\sim \mathcal{I}_{k}^{-1} }\left( \mathcal{I}_{k}\mathcal{I}\right) ^{es}=\sum_{\deg \: f_{k}\overline{g}_0+y{g}_1+ \cdots + y^{q-1} {g}_{q-1}  =d+d_{k}}\left( f_{k}\overline{g}_0+y {g}_1+\cdots +y^{q-1} {g}_{q-1} \right) ^{es}, 
\]
where $\overline{g}_0$ is a polynomial in $x$ and the bar denotes the action of ``dividing by $f_k$''. We examine now, sufficient conditions for the vanishing of the previous sum.
If $q \mu +e_{k} d_{k}<d+d_{k}$, then $d+d_k-e_{k}d_{k}>q \mu >\mu>\frac{l_{q}(es)}{q-1}$ implies that the sum vanishes. So, without loss of generality $d+d_k \leq q\mu +e_{k}d_{k}<N$, and these are the only terms which can give non-zero contributions. As $d+d_{k}<N$, then ${g}_1= \cdots= {g}_{q-1}=0$. Then degree $\overline{g}_0=d+d_{k}-e_{k}d_{k} \leq q \mu$.
Therefore, the $k$-th term of the zeta sum (up to $f_k^{-ss_k}$ factor which can be ignored as already mentioned) is the following:

\begin{eqnarray*}
&&\sum_{d=0}^{\infty }X^{d}\sum_{ \deg \: \mathcal{I}=d \atop \mathcal{I}\sim \mathcal{I}_{k}^{-1} }\left( \mathcal{I}_{k}\mathcal{I}\right) ^{es}= \sum_{d=0}^{\infty }X^{d} \sum_{\deg \: f_{k} \overline{g}_0 =d+d_{k}}\left( f_{k} \overline{g}_0 \right) ^{es}\\
&=&X^{e_{k}d_{k}-d_{k}}f_{k}^{es}+\cdots +X^{^{e_{k}d_{k}-d_{k}}+q\mu }\sum_{a_{\mu -1},\ldots,a_{0}\in \mathbb{F}_{q}}f_{k}^{es}\left( x^{\mu }+a_{\mu -1}x^{\mu -1}+\cdots +a_{0}\right)^{es} \\
&=&f_{k}^{es}X^{e_{k}d_{k}-d_{k}}\left( 1+\cdots +X^{q\mu }\sum_{a_{\mu -1},\ldots,a_{0}\in \mathbb{F}_{q}}\left( x^{\mu }+a_{\mu -1}x^{\mu -1}+\cdots +a_{0}\right) ^{es}\right) \\
&=&f_{k}^{es}X^{e_{k}d_{k}-d_{k}} \left( \sum_{d=0}^{\mu }X^{qd}\sum_{\deg  _{x}\overline{g}_0 =d}\overline{g}_0^{es} \right) \\
&=&f_{k}^{es}X^{e_{k}d_{k}-d_{k}}\zeta _{\mathbb{F}_{q}\left[ x\right] }\left(-es,X^{q}\right).
\end{eqnarray*}
In order to prove the last equality, observe that by theorem \ref{dinesh}, since $\mu>\frac{l_{q}(es)}{q-1}$ implies that
\[
\sum_{d=0}^{\mu }X^{qd}\sum_{\deg _{x}\overline{g}_0=d}\overline{g}_0^{es}=\sum_{d=0}^{\infty }X^{qd}\sum_{\deg _{x}\overline{g}_0=d}\overline{g}_0^{es}=\zeta _{\mathbb{F}_{q}\left[ x\right] }\left(-es,X^{q}\right).
\] 
Now, since the order of vanishing of $\zeta _{\mathbb{F}_{q}\left[ x\right] }\left(-es,X \right)$ is one, we have that, for each $k$, the order of vanishing of the $k$-th term is $q$.\hfill $\Box$

\begin{rem}
\label{exact}
From the proof of theorem \ref{tesismc}, it follows that when $\frac{l_{q}\left( es\right) }{q-1}\leq \mu $ \\ we have
\begin{eqnarray*}
\zeta \left( -es,X\right) =\zeta _{F_{q}\left[ x\right] }\left( -es,X^{q}\right) \left( 1+\sum_{k=1}^{h-1}f_{k}^{\frac{es}{e_{k}}\left( e_{k}-1\right) }X^{e_{k}d_{k}-d_{k}}\right) .
\end{eqnarray*}
Then, the order of vanishing is exactly $q$ when
\begin{eqnarray*}
1+\sum_{k=1}^{h-1}f_{k}^{\frac{es}{e_{k}}\left( e_{k}-1\right) }\neq 0.
\end{eqnarray*}
In particular, when $h=2$ the order of vanishing is exactly $q.$
\end{rem}

In our search for examples of theorem \ref{tesismc} we looked first at $q=2$. We have

\begin{prop}
  Suppose that $q=2$, and $K$ is hyperelliptic. If $K$ has an $r$-gap structure at $\infty$, then $g-1 \leq r \leq g$, where $g$ is the genus of $K$.
\end{prop}

\textbf{Proof: }Assume that $r = g-k$, with $k \geq 2$. Then, since $K$ has an $r$-gap structure, $l\left(2r\infty\right)=l\left(\left(2g-2k\right)\infty\right)=g-k+1$, and $l\left(\left(g+r\right)\infty\right)=l\left(\left(2g-k\right)\infty\right)=g-k+1$. So, between $2g-2k+1 \leq n \leq 2g-k$ all the numbers $n$ are gaps. Observe that, as $k \geq 2$, then there is a gap at an even number, a contradiction.\hfill $\Box$ 

\begin{rem}
Still, we do not have an example of the phenomenon described at theorem \ref{tesismc}. We contacted J. F. Voloch, whose comment was that it is probable that a field satisfying the hypothesis of the theorem exists, but that the method used by him to find the example \ref{voloch} is no longer useful in our case; because, apart from asking for the curve geometric conditions, we are asking also for conditions that are not geometric, for example, the type of curve, the class number, etc. So, he recommended that we do a computational search using the computer program Kant/Kash. But we have not yet been successful.
\end{rem}

\section*{Acknowledgements}

The authors would like to sincerely thank the referee who made valuable stylistic suggestions and also found many errors in earlier versions of this paper. We thank too Gabriel Villa Salvador for his useful comments.

\end{document}